\documentclass[12pt]{article}
\usepackage{theorem}
\usepackage{latexsym}
\usepackage{amssymb}
\usepackage{amsmath}
\usepackage{epic}
\def\Bx{\hfill{$\Box$}}

\def\E{{\cal E}}

\def\A{{\cal A}}
\def\B{{\cal B}}

\newtheorem{thm}{Theorem}[section]
\newtheorem{lemma}[thm]{Lemma}
\newtheorem{cor}[thm]{Corollary}
\newtheorem{defn}{Definition}[section]

\begin{document}

\renewcommand{\theequation}{\arabic{section}.\arabic{equation}}
\thispagestyle{empty}

\vskip 20pt
\begin{center}
{\bf REFINED RESTRICTED PERMUTATIONS}
\vskip 15pt
{\bf Aaron Robertson\footnote{
Homepage:  {\tt http://math.colgate.edu/$\sim$aaron/}}
}\\
{\it Department of Mathematics,}
{\it Colgate University,
Hamilton, NY 13346}\\
{\tt aaron@math.colgate.edu}

{\bf Dan Saracino}\\
{\it Department of Mathematics,}
{\it Colgate University,
Hamilton, NY 13346}\\
{\tt dsaracino@mail.colgate.edu}

{\bf Doron Zeilberger\footnote{Supported
in part by the NSF.  Homepage:  {\tt 
http://www.math.rutgers.edu/$\sim$zeilberg/}
 }
}\\
{\it Department of Mathematics,}
{\it Rutgers University,
New Brunswick, NJ 08903}\\
{\tt zeilberg@math.rutgers.edu}
\end{center}
\vskip 30pt
\hfill {\it Dedication:  In memory of Rodica Simion (1955-2000)}
\vskip 10pt
{\it
This article is dedicated to the memory of Rodica Simion, one of the greatest enumerators of
the 20$^{th}$ century. Both derangements  ([SS2]) and
resticted permutations ([SS]) were very dear to her heart, and we are
sure that she would have appreciated the present surprising connections
between these at-first-sight unrelated concepts. }
\vskip 30pt
\begin{abstract}{\footnotesize \noindent
Define $S_n^k(\alpha)$ to be the set of permutations of $\{1,2,\dots,n\}$
with exactly $k$ fixed points
which avoid the pattern $\alpha \in S_m$.  Let
$s_n^k(\alpha)$ be the size of $S_n^k(\alpha)$.
We investigate $S_n^0(\alpha)$ for all $\alpha \in S_3$ as well
as show that $s_n^k(132)=s_n^k(213)=s_n^k(321)$ 
and $s_n^k(231)=s_n^k(312)$ for all $0 \leq k \leq
n$.  }
\end{abstract}
\vskip 20pt
\section{\large Introduction}

Let $\pi \in S_n$ be a permutation of $\{1,2,\dots,n\}$ written
in one-line notation.
Let $\alpha \in S_m$.
We say that $\pi$ {\it contains the pattern $\alpha$} if there exist
indices $i_1,i_2,\dots,i_m$ such that $\pi_{i_1} \pi_{i_2} \dots
\pi_{i_m}$ is equivalent to $\alpha$, where we define equivalence as
follows. Define $\overline{\pi}_{i_j}=
|\{x:\pi_{i_x} \leq \pi_{i_j}, 1 \leq x \leq m\}|$.  If
$\alpha = \overline{\pi}_{i_1} \overline{\pi}_{i_2} \dots
\overline{\pi}_{i_m}$ then
we say that $\alpha$ and $\pi_{i_1} \pi_{i_2} \dots \pi_{i_m}$ are
equivalent.  For example, if $\tau=124635$ then $\tau$ contains
the pattern $132$ by noting that $\tau_2 \tau_4 \tau_5 = 263$
is equivalent to $132$.  We say that $\pi$ is $\alpha$-{\it avoiding} if
$\pi$ does not contain the pattern $\alpha$. In our above example, $\tau$
is $321$-avoiding.

Define $S_n(\alpha)$, $\alpha \in S_m$, to the
the set of $\alpha$-avoiding permutations in $S_n$.  Let $s_n(\alpha)$
be the size of $S_n(\alpha)$. Knuth ([Knu]) showed that, regardless of the
pattern $\alpha \in S_3$,
$s_n(\alpha)=C_n = \frac{1}{n+1} {2n \choose n}$, the 
$n^{\mathrm{th}}$ Catalan
number.
Bijective results are given in [Kra], [Ric], [SS], and [Wes].

We refine the investigation of $S_n(\alpha)$ in the following
fashion.  Let $\alpha \in S_m$.
Define $S_n^k(\alpha)$ to be the set of $\alpha$-avoiding
permutations of $\{1,2,\dots,n\}$, with exactly $k$ fixed points.
Let
$s_n^k(\alpha)$ be the size of $S_n^k(\alpha)$.
We may write $D_n(\alpha)$ and $d_n(\alpha)$
for $S_n^0(\alpha)$ and $s_n^0(\alpha)$, respectively,
since we are dealing with derangements.

\section{\large Similarity Relations, Catalan Sequences,
and Fine's Sequence}
 
We begin our investigation with similarity relations.
A similarity relation, $R$, is a binary relation
on an ordered set
which is reflexive, symmetric, but not necessarily 
transitive, with the condition that if $iRk$
and $i<j<k$ then $iRj$ and $jRk$.  
Furthermore, we have the following definition
about the structure of a similarity relation.

\begin{defn}
A similarity relation, $R$, is said to have $k$
isolated points if
$k$ is the number of $i \in \{1,2,\dots,n\}$ such that
there does not exists $j \neq i$ with $iRj$.  If $k=0$,
we say that the similarity relation is nonsingular.
We denote by $SR_n(k)$ the set of similarity
relations on $\{1,2,\dots,n\}$ with $k$ isolated points.
\end{defn}

There are two
common structures which can be used to view
similarity relations:  graphs and sequences.
We will be using the sequence interpretation of a
similarity relation as given by Strehl in [Str]:
Let $R$ be a similarity relation on $\{1,2,\dots,n\}$.
Then $R$ corresponds to an integer sequence $r_1r_2\dots r_n$ defined
for any $1 \leq i \leq n$
by $r_i=i-j$, where $j$ is the smallest element
of $\{1,2,\dots,n\}$ such that $iRj$.  Throughout
this paper we assume that similarity relations
are defined on $\{1,2,\dots,n\}$.
To this
end, we make the following definition.

\begin{defn}  The set of similarity relations
(on $\{1,2,\dots,n\}$) is given by
$$
SR_n=\{r_1r_2\cdots r_n \, : r_i \in \mathbb{Z}, \, r_1=0 \,
{\mathrm{and}}
\, 0 \leq r_{i+1} \leq r_i+1 \, \mathrm{for} \, 1 \leq i \leq n-1\}.
$$
\end{defn}

For example $SR_3 = \{000, 001, 010, 011, 012\}$.

It is known ([Rog], [Str]) that for $n \geq 1$, $|SR_n(0)| = F_n$,
where $F_n$ is the $n^{\mathrm{th}}$ Fine number.  The first few
values of Fine's sequence are
$0,1,2,6,18,57,186,622,2120,\dots$, a sequence first discovered in [Fin].

Looking at $SR_n(0)$ a little more closely, we see that 
$s \in SR_n(0)$ if and only if $s$ has no occurrence of
$00$ and does not end with $0$.  Hence, we say that a nonsingular
similarity relation has no {\it double zero}, where we consider
an ending $0$ to be a double zero.

We now state some of the results concerning Fine's sequence as
given in [Rog], [Sha], and [Str].

\begin{thm}  Let $C_n = \frac{1}{n+1} {2n \choose n}$
and $F_n$ be
the $n^{th}$ Catalan number and Fine number,
respectively.  We have the following for $n \geq 2$.
\begin{enumerate}
\item $|SR_n(0)| = F_n$
\item $C_n = 2F_n+F_{n-1}$
\item $F_n = \sum_{1 \leq k \leq n/2} {2n-2k-1 \choose n-1}
- {2n-2k-1 \choose n}$
\item $F_n = \frac12\sum_{i=0}^{n-2} \left(\frac{-1}{2} \right)^i C_{n-i}$
\end{enumerate}
\end{thm}

As we can see, the Catalan and Fine numbers are related.
This becomes more evident in light of the following
definition.

\begin{defn}   The set of Catalan sequences
of length $n$ is given by
$$
Cat(n)=\{ c_1c_2 \cdots  c_n \, :\,c_i \in \mathbb{Z}, \,
 1 \leq c_1 \leq c_2 \leq \dots \leq c_n  , \,
{\mathrm{and}} \, c_i \leq i \,\, {\mathrm{for}} \, 1 \leq i \leq n \}
$$

\end{defn}

For example:
$Cat(3)=\{111,112,113,122,123\}$. It is well-known and
easy to see that the cardinality of $Cat(n)$ is the
$n^{\mathrm{th}}$ Catalan number. It is
also well-known and easy to see that the generating function for the
Catalan numbers, $\psi(t)=\sum_{n=0}^{\infty} C_n t^n$, satisfies
the quadratic equation $\psi(t)=1+t\psi^2(t)$, and hence that
$\psi(t)=\frac{1-\sqrt{1-4t}}{2t}$.

\section{\large 321-Avoiding Derangements and Dyck Paths}

The aim of this section is to show that the $321$-avoiding
derangements are enumerated by Fine's sequence.  We will
investigate two bijections, the main one
due to Krattenthaler [Kra].  First,
we must introduce a few definitions.

\begin{defn} We say that a permutation $\pi \in S_n$ is
a backward derangement if $\pi_{n+1-i} \neq i$ for
all $1 \leq i \leq n$.
\end{defn}

Consequently, a $123$-avoiding backward derangement
when read from right to left is a $321$-avoiding derangement.

We will be using a bijection
due to Krattenthaler [Kra] from $S_n(123)$ to the
set of Dyck paths of length $2n$, so for completeness
we define a Dyck path.

\begin{defn} A Dyck path is a path in $\mathbb{R}^2$ from
$(0,0)$ to $(2n,0)$ consisting of a sequence of steps
of length $\sqrt{2}$ and slope $\pm 1$.  We denote
these two types of steps by $(1,1)$ and $(1,-1)$,
called up-steps and down-steps, respectively.  We say that
the length of such a Dyck path is $2n$ (its horizontal
length) and denote the set of Dyck paths of length
$2n$ by $Dyck(2n)$.
\end{defn}

We also have the following definition about certain Dyck
paths.

\begin{defn}  We say that a Dyck path contains a hill if it
has a peak at height $1$.  We say that a Dyck path is
hill-free is it contains no hill.  We denote the set of
hill-free Dyck paths of length $2n$ by $Dyck^{hf}(2n)$.
\end{defn}

We now describe a natural bijection from $SR_n$ to  Dyck
paths of length $2n$.  Let $s=s_1s_2\dots s_n \in SR_n$.
This bijection is very similar to one given by Krattenthaler
[Kra] from $S_n(132)$ to $Dyck(2n)$.

Each $s_i$, $1 \leq i \leq n$, corresponds to the starting height of
an up-step.  Proceeding from $s_i$ to $s_{i+1}$, if $s_{i+1} > s_i$
then we continue with up-steps.  If $s_{i+1} \leq s_i$ we
append $s_i-s_{i+1}+1$ down-steps followed by a single up-step.
This assures us that $s_{i+1}$ corresponds to a starting height
of $s_{i+1}$
for an up-step.  After $s_n$ we use as many down-steps as
necessary to end at $(2n,0)$. 
The inverse bijection is obvious.  An example is in order.

Let $s=0120121 \in SR_7(0)$.  This nonsingular similarity
relation corresponds to the following Dyck path, with
the entries of $s$ marked on the Dyck path.

\hskip -33pt
\setlength{\unitlength}{1mm}
\begin{picture}(0,45)(-20,0)
\linethickness{.5mm}
%axes
\put(0,0){\line(0,1){40}}
\put(0,0){\line(1,0){150}}
\linethickness{.1mm}
%thin lines
\multiput(0,0)(0,10){4}{\line(1,0){150}}
\multiput(0,0)(10,0){15}{\line(0,1){40}}
%Labels
\put(-5,-1){0}
\put(-5,9){1}
\put(-5,19){2}
\put(-5,29){3}
\put(-1,-5){0}
\put(9,-5){1}
\put(19,-5){2}
\put(29,-5){3}
\put(39,-5){4}
\put(49,-5){5}
\put(59,-5){6}
\put(69,-5){7}
\put(79,-5){8}
\put(89,-5){9}
\put(99,-5){10}
\put(109,-5){11}
\put(119,-5){12}
\put(129,-5){13}
\put(139,-5){14}
%vertices
\put(-1,-1){$\bullet$}
\put(9,9){$\bullet$}
\put(19,19){$\bullet$}
\put(29,29){$\bullet$}
\put(39,19){$\bullet$}
\put(49,9){$\bullet$}
\put(59,-1){$\bullet$}
\put(69,9){$\bullet$}
\put(79,19){$\bullet$}
\put(89,29){$\bullet$}
\put(99,19){$\bullet$}
\put(109,9){$\bullet$}
\put(119,19){$\bullet$}
\put(129,9){$\bullet$}
\put(139,-1){$\bullet$}
\linethickness{.5mm}
\put(0,0){\line(1,1){30}}
\put(30,30){\line(1,-1){30}}
\put(60,0){\line(1,1){30}}
\put(90,30){\line(1,-1){20}}
\put(110,10){\line(1,1){10}}
\put(120,20){\line(1,-1){20}}
%labels
\put(3,5){{\bf 0}}
\put(13,15){{\bf 1}}
\put(23,25){{\bf 2}}
\put(63,5){{\bf 0}}
\put(73,15){{\bf 1}}
\put(83,25){{\bf 2}}
\put(113,15){{\bf 1}}
\end{picture}
\vskip 25pt
\small
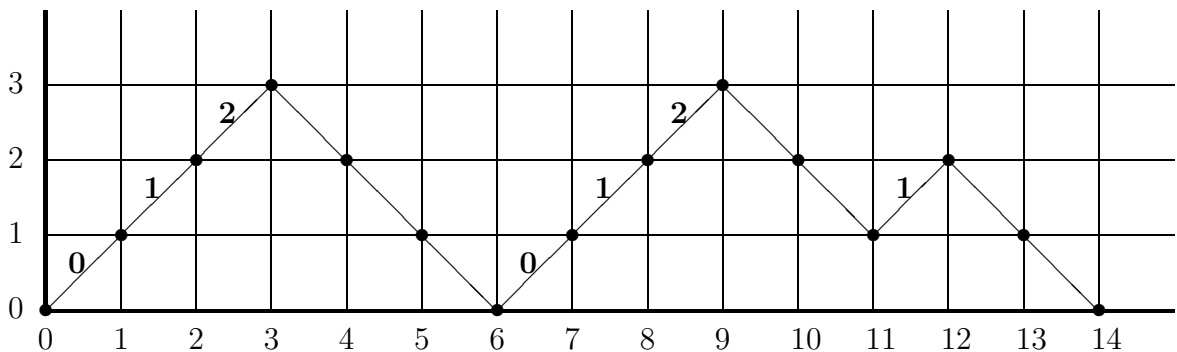
\centerline{{\bf Dyck path corresponding
to $\mathbf{0120121 \in SR_7(0)}$ and to $\mathbf{6573142 \in S_7(123)}$}}
\normalsize

\vskip 20pt

Using the above bijection and the fact that $|Dyck(2n)|=C_n$,
we easily obtain $|SR_n|=C_n$
(which was shown in [Rog] and [Str]).  Furthermore, we get
the following theorem ($SR_n(0)=F_n$ was shown in
[Rog] and [Str] while  $|Dyck^{hf}(2n)|=F_n$
was shown in [Deu] ).  

\begin{thm} For $n \geq 1$, $|SR_n(0)|=|Dyck^{hf}(2n)|=F_n$, 
where $F_n$ is the $n^{{th}}$ Fine number
\end{thm}

{\bf Proof.}  Clearly, a hill occurs in a Dyck path
if and only if the corresponding similarity relation contains
a double zero.
\Bx

We are now ready to use Krattenthaler's bijection ([Kra]),
described as follows.
Let $\pi=\pi_1 \pi_2 \cdots \pi_n \in S_n(123)$.  Determine the
right-to-left maxima of $\pi$, i.e. $m=\pi_i$ is a 
right-to-left maximum if $m>\pi_j$ for all $j>i$.  For example, the
right-to-left maxima of
$65\mathbf{7}31\mathbf{4}\mathbf{2}$ are in bold type.

Let $\pi$ have right-to-left maxima $m_1<m_2<\dots<m_s$, so that
we may write
$$
\pi=w_sm_sw_{s-1}m_{s-1}\cdots w_1m_1.
$$
So, for example, in
$6573142$ we have $w_3=65$, $m_3=7$, $w_2=31$, $m_2=4$,
$w_1 = \emptyset$, $m_1=2$.

We now generate a Dyck path from $(2n,0)$ to $(0,0)$
(generate it backwards) using {\it backward} up-steps $((-1,1))$
and {\it backward} down-steps $((-1,-1))$.

Read $\pi$ from right to left.  For each $m_i$ do $m_i-m_{i-1}$
up-steps (where we define $m_0=0)$.  For each $w_i$ do $|w_i|+1$
down-steps.  So, for example, $\pi=6573142$ generates
the Dyck path shown above.

We now state the main theorem of this section.

\begin{thm}
Let $\pi \in S_n(123)$ and let $D_\pi$ be the associated Dyck
path provided by Krattenthaler's bijection.  Then
$\pi$ is a backward derangement if and only if $D_\pi$
is hill-free.
\end{thm}

\noindent
{\bf Proof.}
We first demonstrate the ``only-if" direction.
Assume, for a contradiction, that $D_\pi$ has a hill.
Notice that the hill cannot be at either end of the Dyck
path for otherwise $\pi_1=n$ or $\pi_n=1$, both contradicting
the fact that $\pi$ is a backward derangement.  Hence,
any hill must be an ``interior" hill.  Let $\pi=\pi(1)\pi(2)$
and $D_\pi
= D_{\pi(1)} D_{\pi(2)}$ where $D_{\pi(1)}$ contains
a hill on the right end and $D_{\pi(2)}$ contains
no hill.  Let $D_{\pi(2)}$ consist of $2y$ steps.  
By construction, $\pi(2)$ contains the elements $1,2,\dots,y$.
Thus, in order for $D_{\pi(1)}$ to end with a hill, we
must have $\pi_{n-y}(1)=y+1$, i.e. $\pi_{n+1-(y+1)}=y+1$, contradicting
the fact that $\pi$ is a backward derangement.

We now give the ``if" direction; we prove the contrapositive.
Assume, for a contradiction, that $D_\pi$
has no hill.
Let $x$ be the smallest integer such that
$\pi_{n+1-x}=x$.  Write $\pi=\pi(1)x\pi(2)$.
 From the above argument, we see
that $\pi(2)$ cannot consist of the elements
$1,2,\dots,x-1$ only, for otherwise we would
have a hill.
Hence, there exists $y \in \pi(2)$ such
that $y>x$.
Consequently, there exists $w \in \pi(1)$
such that $w<x$.  But then $wxy$ is a $123$-pattern,
a contradiction.
\Bx

Coupling this theorem with Theorem 3.1,
 we immediately obtain the
following corollary.

\begin{cor} For $n \geq 1$, $d_n(321)=F_n$, where
$F_n$ is the $n^{th}$ Fine number.
\end{cor}

We can now investigate the relationship between
$S_n^k(321)$ and similarity relations.
To this end, we have the following theorem.

\begin{thm} For $n \geq 0$, $s_n^k(321)=|SR_n(k)|$.
\end{thm}

\noindent
{\bf Proof.}  We prove this by induction on $k$.
The case $k=0$ is given by Corollary 3.3.  Hence,
we assume that there exists a bijection
$\gamma_n^k: S_n^k(321) \rightarrow SR_n(k)$.

Let $\pi \in S_n^{k+1} (321)$ and let $f$ be the
smallest fixed point.  Write $\pi=\pi(1)f\pi(2)$
so that $\pi(2)$ contains $k$ fixed points.
Note that since $\pi$ must be $321$-avoiding we must
have $\pi(1) \in S_{f-1}^0(321)$ on the elements $1,2,\dots,f-1$
and $\pi(2) \in S_{n-k}^k (321)$ on the elements
$f+1,f+2,\dots,n$.
Let $\gamma_{f-1}^0(\pi(1))=t \in SR_f(0)$ and 
$\gamma_{n-f}^k(\pi(2))=r\in SR_{n-f}(k)$.
Then define $\gamma_n^{k+1} (\pi) = t0r \in SR_n(k+1)$.
To show this is a bijection, it is enough to
give the inverse.  This is obtained by noting that
the position of the first $0$ in the first occurrence of
a double zero in an element of $SR_n(k+1)$ corresponds to the minimal
fixed point.
\Bx

\section{\large Weighted-Counting of 321-Avoiding Permutations}

In order to achieve one of our goals we must enumerate
$S_n^k(321)$ via another approach.  Before delving into
our approach, we make the following definition.

\begin{defn}  Let $S$ be a finite set with each
element $s \in S$ having a unique characteristic from
$C = \{c_1,c_2,\dots,c_k\}$, written as $char(s)$.  The weight-enumerator
of $S$ with respect to $weight(s)=x^{char(s)}$
is given by
$$
\sum_{i=1}^k s_i x^{c_i},
$$
where $s_i = \vert \{s \in S : char(s)=c_i\} \vert$.
\end{defn}

Applying this to our situation,
let $\A_n=S_n(321)$  and let
$A_n(x)$ be its weight-enumerator with respect to
$weight(\pi)=x^{f(\pi)}$, where $f(\pi)$ is the number
of fixed points of $\pi \in S_n(321)$.

Recall that $Cat(n)$ is the set of Catalan sequences,
defined in Section 2.
We will define a bijection $T:\A(n) \rightarrow Cat(n)$ as follows.
If $n=1$, then $T(1)=1$.
For $n \geq 2$ we define $c=T(\pi)$ recursively as follows.

Given a $321$-avoiding permutation $\pi$ of length $n$,
let $i$ be the place where $n$ is (i.e. $\pi_i=n$).
If $\pi_n=n-1$ then let $\pi'$ be 
$\pi_1, \dots, \pi_{i-1},n-1, \pi_{i+1}, \dots, \pi_{n-2}$, 
otherwise let $\pi'$ be $\pi$ with $n$ removed, i.e.
$\pi'=\pi_1, \dots, \pi_{i-1},\pi_{i+1}, \dots, \pi_n$.
Then $c=T(\pi)$ is defined to be $c'=T(\pi')$ with $i$ appended
at the end.

The inverse bijection  $S:Cat(n) \rightarrow \A(n)$  is defined
as follows. If $n=1$ then
$S(1)=1$. If $n \geq 2$ and $c=c_1 \cdots c_n$ is a Catalan sequence
then $\pi=S(c)$ is defined recursively as follows.
Let $c_n=i$, $c_{n-1}=j$, and
let $c'=c_1 \cdots c_{n-2}c_{n-1}$ be $c$ with its
last component removed.
Let $\pi'=S(c')$. If $i\leq  j$ then
let $\pi$ be the permutation obtained from $\pi'$ by
changing the $n-1$ into $n$ and appending $n-1$ to the end,
while if $j<i$ then let $\pi$ be the permutation obtained
from $\pi'$ by inserting $n$ at the $i^{\mathrm{th}}$ place
i.e. $\pi=\pi'_1, \dots, \pi'_{i-1},n,\pi'_{i}, \dots, \pi'_{n-1}$.

It is easy to prove, by induction on
$n$, that $TS$ and $ST$ are identity mappings, and hence
that $T$ is indeed a bijection. It is also easy to see
that if  $c_1\cdots c_n=T(\pi)$ then 
for $1 \leq i <n$, $\pi_i=i$ if and only if $c_{i}=i$ and $c_{i+1}=i+1$,
and $\pi_n=n$ if and only if $c_n=n$. Hence $A_n(x)$ is equal to
the weight-enumerator of Catalan sequences with
respect to $weight(c)=
x^{g(c)}$, where $g(c)$ is the number of $i$'s such that
$c_{i}=i$ and $c_{i+1}=i+1$ ($i<n$)
plus 1 if $c_n=n$.

Given a Catalan sequence $c$, let
$$ 
D(c)=\{ i \, : 1 \leq i< n, c_i=i , c_{i+1}=i+1 \, \mathrm{or} \,
i=n \, \mathrm{and} \, c_n=n \}.
$$
For example, $D(11345558)=\{3,4,8\}$, 
$D(11111111)=\emptyset$,  and
$D(112346)=\{6\}$. Note the weight of a Catalan sequence $c$ is
$x^{\vert D(c)\vert }$.

In order to weight-enumerate the set of Catalan sequences
it would be easier to use the inclusion-exclusion
philosophy and consider the larger sets of {\it marked}
Catalan sequences, which are the sets of pairs
$(c,S)$ with $S \subset D(c)$ and weight defined
by $weight(c,S)=(x-1)^{\vert S \vert }$.  For example
$weight(11345558, \{3,8\})=(x-1)^2$ and
$weight(11345558, \{\})=(x-1)^0=1$.
Since $x^{\vert D(c)\vert}=((x-1)+1)^{\vert D(c)\vert}$=
$\sum_{S \subset D(c)} (x-1)^{\vert S \vert}$, it follows
that $A_n(x)$ is the weight-enumerator of marked Catalan
sequences.

We now derive a recurrence. Given a marked Catalan sequence
$(c_1 \cdots c_n,S)$, if $1 \in S$ (i.e. $1$ is marked)
then we can get a smaller marked word by deleting
$c_1=1$ and diminishing all indices and elements of $S$ by $1$. The
weight-enumerator of this case is $(x-1)A_{n-1}(x)$.

If $1 \not \in S$ (i.e. $1$ is not marked), let
$i$ be the smallest $i>1$ such that $c_i=i$, if it exists.
Then $c_2 \cdots c_{i-1}$ is a run-of-the-mill
Catalan sequence of length $i-2$, 
while 
$$
\left(c_i-(i-1) \,c_{i+1}-(i-1) \, \cdots\, c_n-(i-1),S-(i-1)\right)
$$
(by $S-(i-1)$ we mean the set $S$ with all of its elements reduced by
$i-1$) is a marked Catalan sequence of length $n-i+1$.  The
weight-enumerator of this, for a given $i$, is
$C_{i-2}A_{n-i+1}(x)$ $(2 \leq i \leq n)$. Finally, if
$c_i<i$ for all $2 \leq i \leq n$ then 
$c_2 \cdots c_n$ is a typical Catalan sequence
of length $n-1$. Hence we get that
$$
A_n(x)=(x-1)A_{n-1}(x)+\left ( 
\sum_{i=2}^{n} C_{i-2}A_{n-i+1}(x) \right )+C_{n-1},
$$
which can be rewritten as
\begin{equation}
A_n(x)=(x-1)A_{n-1}(x)+\sum_{i=1}^{n} C_{i-1}A_{n-i}(x).
\end{equation}

Introducing the generating function
$$
\phi(x,t)=\sum_{n=0}^{\infty} A_n(x)t^n,
$$
and recalling the generating function for the Catalan numbers
$\psi(t)=\frac{1-\sqrt{1-4t}}{2t}$,
$(4.1)$ translates to
$$
\phi(x,t)=1+(x-1)t\phi(x,t)+t\phi(x,t)\psi(t).
$$
Solving for $\phi(x,t)$ yields the explicit expression
$$
\phi(x,t)=\frac{2}{1-2(x-1)t+\sqrt{1-4t}} 
$$
which, in turn, implies (by multiplying the top and the bottom by
$1-2(x-1)t-\sqrt{1-4t}$) that
$$
\phi(x,t)=\frac{1-x+\psi(t)}{2-x+t(x-1)^2}.
$$

To sum up we have the following theorem.

\begin{thm} Let $a_n(x)$ be the coefficient of $t^n$ in
the Maclaurin expansion 
with respect to $t$ of
$$
\frac{1-x+\psi(t)}{2-x+t(x-1)^2}.
$$
Then $A_n(x)=a_n(x)$.
\end{thm}

We now investigate an interesting property of
$a_n(x)$ that will be used in the next section.
Expanding $a_n(x)$ in powers of $(x-1)$ and
using (4.1), we may write

\begin{equation}
a_n(x) = (x-1)^n+\sum_{i=0}^{n-1}\sum_{j=1}^{n-i} C_{j-1}a_{n-i-j}(x)
\cdot (x-1)^i.
\end{equation}

Using (4.2) we prove the following 
crucial lemma.

\begin{lemma}  Let $a_n(x) = \sum_{k=0}^n b(n,k)(x-1)^k$
and define $b(x,y)=0$ if $y>x, x<0$, or $y<0$.
Then, for $n \ge 1$ and $k \neq -1$, 
\begin{equation}
b(n,k)=b(n,k+1)+b(n-1,k-1)
\end{equation}
where $b(0,0)=1$.

\end{lemma}

\noindent
{\bf Proof.} Note that to render
(4.3) valid when $n=k=0$, we would have to add a correction
term of $1$ to the right side, and to render (4.3) valid
when $n \geq 0$ and $k=-1$ we would have to add
$b(n,0)=C_n$ to the left side.

We now prove (4.3) for $n \geq 1$ and $k\neq -1$ by induction on
$n$. Since $a_0(x)=1$ and $a_1(x)=x$, the case $n=1$ is true.
We now assume that
$b(n-i,k)=b(n-i,k+1)+b(n-i-1,k-1)$ for all $1 \leq i \leq n-1$
and $k \neq -1$.

To establish (4.3) for $n$ and all $k \neq -1$ we first note that
this is clear if $k<-1$ or $k>n$.  It is also clear for $k=n$
since $b(n,n)=b(n-1,n-1)=1$ by (4.2).  To deal with
$0 \leq k \leq n-1$ we use the fact
that (4.2) yields
\begin{equation}
b(n,k) = \sum_{i=0}^{n-1} \sum_{j=1}^{n-i} 
C_{j-1}b(n-i-j,k-i)
\end{equation}
for $0 \leq k \leq n-1$.

In particular, for $k=n-1$ this says that $b(n,n-1)=n$,
so $b(n,n-1)=b(n,n)+b(n-1,n-2)$ and (4.3) holds for $k=n-1$.

For $0 \leq k \leq n-2$ we take (4.4) and apply
(4.3) to each term on the right side, adding a correction term of
$C_{n-k-1}$ when $i=k$ and $j=n-k$
(note that this case occurs since $k \neq -1$), and
correction term $-\sum_{i=1}^{n-k-1} C_{i-1}C_{n-k-i-1}$
when $i=k+1$ (note that this case occurs since $k \neq n-1$).
Since the correction terms cancel, we have
$$
b(n,k) = \sum_{i=0}^{n-1} \sum_{j=1}^{n-i} C_{j-1}b(n-i-j,k-i+1) +
\sum_{i=0}^{n-1} \sum_{j=1}^{n-i} C_{j-1}b(n-i-j-1,k-i-1).
$$
Since $k+1 \leq n-1$, the first term on the right is $b(n,k+1)$
by (4.4).  By (4.4) again, the second term on the right is $b(n-1,k-1)$
except for the terms in the double sum where $i=n-1$
or $j=n-i$.  But these terms are all $0$ because they involve 
values of $b(x,y)$ with $x<0$.  This concludes the proof.
\Bx

In [Str], Strehl defines $a(n,k)$, for $n \geq 1$ and $1 \leq k \leq n$,
as the number of similarity relations on $\{1,2,\dots,n\}$
that have $k$ zeros.  He proves that
$a(n,k)=a(n,k+1)+a(n-1,k-1)$ for $n \geq 2$, with
$a(n,1)=a(n,2)=C_{n-1}$.  It follows from Lemma 4.2
that $b(n,k)=a(n+1,k-1)$.

Using Lemma 4.2 we see that the beginning of the $b(n,k)$
table is as follows, where $n=0,1,\dots$ corresponds to row $n$
and $k=0,1,2,\dots$ is the $k^{\mathrm{th}}$ term in from the left. 
This generates the so-called Catalan Triangle ([Slo]).  
\vskip -15pt
\small
$$
\begin{array}{ccccccccccccc}
&&&&&&1&&&&&\\
&&&&&1&&1&&&&\\
&&&&2&&2&&1&&&\\
&&&5&&5&&3&&1&&\\
&&14&&14&&9&&4&&1&\\
&42&&42&&28&&14&&5&&1\\
132&&132&&90&&48&&20&&6&&1\\
\end{array}
$$
\vskip 10pt
\centerline{\bf Beginning of the Catalan triangle; values of
$\mathbf{b(n,k)}$}
\vskip 20pt
\normalsize
We also use Lemma 4.2 to prove the next lemma, which
will be a crucial step in
one of our main proofs.

\begin{lemma} Let $a_n(x) = \sum_{k=0}^n b(n,k)(x-1)^k$.
Then, for all $n,k \geq 1$,
$$
b(n,k)=\sum_{i=1}^n C_{i-1}b(n-i,k-1)
$$ 
where $b(n,0)=C_n$, 
and $b(x,y)=0$ if $y>x$.
\end{lemma}

\noindent
{\bf Proof.}  We use double induction; forward on $n$
and backward on $k$.  We start with induction on $n$.
From the table above, this clearly holds for $n=2$ and $1 \leq k 
\leq 2$.  Hence, we assume that it holds for $n-1$
and $1 \leq k \leq n-1$ to show that it holds
for $n$ and $1 \leq k \leq n$.  We now perform backward
induction on $k$.  The base case $b(n,n)$
holds since $b(n,n)=b(n-1,n-1)=1$.
Hence, we assume that the case $b(n,k+1)$ holds to show that
the case $b(n,k)$ holds.

We must do the case $k=1$ seperately.  This holds by
the identity $C_n = \sum_{i=1}^n C_{i-1}C_{n-i}$
and the fact that $b(n,1)=b(n,0)=C_n$.

Now, from Lemma 4.2, we have $b(n,k) = b(n,k+1)+b(n-1,k-1)$
rendering the induction straightforward for $k\geq 2$.
\Bx

\section{\large Weighted-Counting of 132-Avoiding Permutations}
\setcounter{equation}{0}

Recall that $\A_n=S_n(321)$
and 
$A_n(x)$ is its weight-enumerator with respect to
$weight(\pi)=x^{f(\pi)}$, where $f(\pi)$ is the number
of fixed points of $\pi \in S_n(321)$.
Let $\B_n=S_n(132)$  and let
$B_n(x)$ be its weight-enumerator with respect to the same weight.
Our goal is to show that 
for $n \geq 0$, $A_n(x)=B_n(x)$.

We have that $\B_n$ is the set of $132$-avoiding
bijections on $\{1,2, \dots, \, n \}$.
However, in order to weight-enumerate them, we must
consider, more generally, the set $\E_{n,r}$ of bijections
$$
\pi:\{1,2, \dots, \, n \} \rightarrow \{r+1,r+2, \dots, r+n \} 
$$
that avoid $132$. Let $E_{n,r}(x)$ be the weight-enumerator of
these bijections, with respect to the same weight. It
is easy to see that
$E_{n,-r}(x)=E_{n,r}(x)$, and that $E_{n,r}(x)=C_n$ if
$r \geq n$, since in this situation no fixed point can occur.
Of course $B_n(x)=E_{n,0}(x)$.

We will now establish a recurrence for the $E_{n,r}(x)$.
Consider $\pi \in \E_{n,r}$.
Let $i$ be the location of $n+r$, i.e.
the index $i$ for which $\pi_i=n+r$.
Now the set of entries
before the $i^{\mathrm{th}}$ place must consist of the
$i-1$ largest elements of the range, and the entries
after the $i^{\mathrm{th}}$ place must consist of the $n-i$ smallest
elements, since otherwise a delinquent $132$ will be formed.
Hence, every member $\pi \in \E_{n,r}$ gives rise
to a pair $(\pi', \pi'')$ where
$\pi'\in \E_{i-1,n-i+r}$ and
$\pi''\in \E_{n-i,r-i}$ . Furthermore, the weight of $\pi$
is the product of the weights of $\pi'$ and $\pi''$, except
when $r=0$ and $i=n$ in which case we removed a fixed point,
and we have an extra factor of $x$.
Let $I_S$ be the characteristic function, i.e.
$I_S=1$ if $S$ is true, and $0$ otherwise.
The above argument gives the following non-linear recurrence.
$$
E_{n,r}(x)=
\sum_{i=1}^{n}
E_{i-1,n-i+r}(x)E_{n-i,r-i}(x)
+(x-1) E_{n-1,0}(x) I_{\{r=0\}}.
$$

Changing $E_{n-i,r-i}(x)$ to $E_{n-i,i-r}(x)$ when $i>r$ and using the
fact that $E_{i-1,n-i+r}(x)=C_{i-1}$ when $n-i+r \geq i-1$,
i.e. $i \leq \lfloor \frac{n+r+1}{2} \rfloor$, and
$E_{n-i,i-r}(x)=C_{n-i}$ when $i-r \geq n-i$, i.e. $i \geq \lceil
\frac{n+r}{2} \rceil$, we get a simplified recurrence, which is 
linear in $E_{n,r}$:

\begin{equation*}
\begin{split}
E_{n,r}(x)&=
\sum_{i=1}^{r} C_{i-1}E_{n-i,r-i}(x)
+\sum_{i=r+1}^{\lfloor \frac{n+r+1}{2} \rfloor}  C_{i-1}E_{n-i,i-r}(x)\\
\\ &+\sum_{i=\lfloor \frac{n+r+1}{2} \rfloor+1}^{n}
E_{i-1,n-i+r}(x)C_{n-i} +(x-1) E_{n-1,0}(x) I_{\{r=0\}}.
\end{split}
\end{equation*}

We are now able to prove the following theorem.

\begin{thm} Let $a_n(x)$
be the weight-enumerator of $321$-avoiding permutations
with respect to  $weight(\pi)=x^{f(\pi)}$, where $f(\pi)$ is
the number of fixed points in $\pi \in S_n(321)$.
Let $E_{n,r}(x)$
be the weight-enumerator of $132$-avoiding bijections
from $\{1,2,\dots,n\}$ to $\{r+1,r+2,\dots,r+n\}$
with respect to the same weight.

Then, for $r \geq 0$,
\begin{equation}
E_{n,r}(x)=a_n(x)+(1-x) \sum_{i=1}^{r} C_{i-1} a_{n-i}(x),
\end{equation}
where we define $a_j(x) =0$ for $j<0$.
\end{thm}

\noindent
{\bf Proof.} 
Denote by $e_{n,r}(x)$ the right side of (5.1).  We must show
for all $n$ (keeping $r$ as a parameter) that
\begin{equation*}
\begin{split}
e_{n,r}(x)&=\sum_{i=1}^{r} C_{i-1}e_{n-i,r-i}(x)
+\sum_{i=r+1}^{\lfloor \frac{n+r+1}{2} \rfloor}  C_{i-1}e_{n-i,i-r}(x) \\
&+
\sum_{i=\lfloor \frac{n+r+1}{2} \rfloor +1}^{n} e_{i-1,n-i+r}(x)C_{n-i}
+(x-1) e_{n-1,0}(x) I_{\{r=0\}}.
\end{split}
\end{equation*}
In other words, we have to prove
\begin{equation}
\begin{split}
a_n(x)+(1-x) \sum_{i=1}^{r} C_{i-1} a_{n-i}(x)
&=
\sum_{i=1}^{r} C_{i-1}
\left ( a_{n-i}(x)+(1-x) \sum_{j=1}^{r-i} C_{j-1} a_{n-i-j}(x) \right )\\
&+\sum_{i=r+1}^{\lfloor \frac{n+r+1}{2} \rfloor}  C_{i-1} 
\left ( a_{n-i}(x)+(1-x) \sum_{j=1}^{i-r} C_{j-1} a_{n-i-j}(x) \right )\\
&+ \sum_{i=\lfloor \frac{n+r+1}{2} \rfloor +1}^{n} 
C_{n-i} 
\left ( a_{i-1}(x)+(1-x) \sum_{j=1}^{n-i+r} C_{j-1} a_{i-j-1}(x) \right
)\\ 
&+(x-1)a_{n-1}(x) I_{\{r=0\}},
\end{split}
\end{equation}
which has been checked by the Maple package {\tt AARON} for $n \leq 50$.
{\tt AARON} was
written by the third author (with
additions by the first author) and
is available at each of these author's website
(given on the first page).

Let the left side of (5.2) be $l_n(x)$
and the right side be $r_n(x)$.
Write $$
l_n(x) = \sum_{k=0}^n s(n,k) (x-1)^k \,\, {\mathrm{and}} \,\,
r_n(x) =  \sum_{k=0}^n t(n,k) (x-1)^k.
$$
If we can show that for any $n$, $s(n,k)=t(n,k)$
for all $0 \leq k \leq n$, then we will be done.

We must take care of the case $k=0$ seperately.  This case  holds
since $b(n,0)=C_n$ and $C_n =\sum_{i=1}^n C_{i-1}C_{n-i}$.

We now use double induction on $n \geq 2$ and $k \geq 0$; forward
induction on $n$ and backward induction on $k$.

We start by inducting on $n$.
Since Maple has given us the base case, i.e.
$s(2,k)=t(2,k)$ for $0 \leq k \leq 2$, we may assume
that $s(n-1,k)=t(n-1,k)$ for all $0 \leq k \leq n-1$.
We must show that $s(n,k)=t(n,k)$ for all $0 \leq k \leq n$.

We proceed via backward induction on $k$.
For our base case we must show that $s(n,n)=t(n,n)$.
Gathering the $(x-1)^n$ terms in (5.2) we have
$(x-1)^n - (x-1)^n I_{\{r \neq 0\}}$ on the left side
of (5.2), and $(x-1)^n I_{\{r =0\}}$ on the right side
of (5.2).  Since $(x-1)^n = (x-1)^n I_{\{r \neq 0\}}
+ (x-1)^n I_{\{r =0\}}$, we have $s(n,n)=t(n,n)$

We now assume that $s(n,k+1)=t(n,k+1)$ to show
that $s(n,k)=t(n,k)$ for $k \geq 1$.  Letting $a_n(x) = \sum_{k=0}^n
b(n,k)(x-1)^k$, with $b(0,0)=1$ and $b(x,y) =0$ if $y>x$, $x<0$, or $y<0$,
we must show that
\begin{equation}
\begin{split}
b(n,k)- \sum_{i=1}^{r} C_{i-1} b(n-i,k-1)
&=
\sum_{i=1}^{r} C_{i-1}
\left ( b(n-i,k)-\sum_{j=1}^{r-i} C_{j-1} b(n-i-j,k-1) \right )\\
&+\sum_{i=r+1}^{\lfloor \frac{n+r+1}{2} \rfloor}  C_{i-1} 
\left ( b(n-i,k)- \sum_{j=1}^{i-r} C_{j-1} b(n-i-j,k-1) \right)\\
&+ \!\!\!\! \sum_{i=\lfloor \frac{n+r+1}{2} \rfloor +1}^{n} 
\!\!\!\!\!\! C_{n-i} 
\left ( b(i-1,k)- \!\!\sum_{j=1}^{n-i+r} C_{j-1} b(i-j-1,k-1) \right
)\\ 
&+b(n-1,k-1) I_{\{r=0\}}.
\end{split}
\end{equation}

The cases $r \geq n$ are straightforward and are left to the
reader.  Below, we assume that $r<n$.

Using Lemma 4.2, and the inductive hypothesis for $b(n,k+1)$
we must show that 
\begin{equation*}
\begin{split}
b(n-1,k-1)- \sum_{i=1}^{r} C_{i-1} b(n-i-1,k-2)
&=
\sum_{i=1}^{r} C_{i-1}
 b(n-i-1,k-1)\\
&-\sum_{i=1}^{r}\sum_{j=1}^{r-i} C_{i-1}C_{j-1} b(n-i-j-1,k-2)\\
&+\sum_{i=r+1}^{\lfloor \frac{n+r+1}{2} \rfloor}  C_{i-1} 
 b(n-i-1,k-1)\\
&- \sum_{i=r+1}^{\lfloor \frac{n+r+1}{2} \rfloor}\sum_{j=1}^{i-r} 
C_{i-1}C_{j-1}
b(n-i-j-1,k-2) \\ 
&+ \sum_{i=\lfloor \frac{n+r+1}{2} \rfloor +1}^{n} 
C_{n-i} 
b(i-2,k-1)\\
&-\sum_{i=\lfloor \frac{n+r+1}{2} \rfloor +1}^{n}\sum_{j=1}^{n-i+r}
C_{n-i}C_{j-1} b(i-j-2,k-2)\\  
&-C_{\lfloor\frac{n+r+1}{2}
\rfloor -1} C_{n - \lfloor
\frac{n+r+1}{2} \rfloor-1} I_{\{k=1\}}\\
&+b(n-2,k-2) I_{\{r=0\}},
\end{split}
\end{equation*}
where the second to last term occurs for $i = \lfloor
\frac{n+r+1}{2} \rfloor $ and $j=n-i$ in the
second double sum of (5.3)
(as this gives $b(0,0)$ which is
equal to $1$ and not $0$ as would be given by the
recurrence in Lemma 4.2).

We must consider two cases:  $n+r$ even and $n+r$ odd.
Assume that $n+r$ is even; the case where $n+r$ is odd is
similiar.  For $n+r$ even we have
$\lfloor \frac{n+r+1}{2} \rfloor = \lfloor \frac{n+r}{2} \rfloor$.
Using the inductive hypothesis for $b(n-1,k-1)$ we are reduced to showing
that
$$
C_{n - \lfloor \frac{n+r}{2} \rfloor-1} b(\lfloor \frac{n+r}{2} \rfloor
-1,k-1)=
C_{\lfloor\frac{n+r}{2}
\rfloor -1} C_{n - \lfloor
\frac{n+r}{2} \rfloor-1}
$$
if $k=1$, which holds
since $b(n,0)=C_n$, or
$$
C_{n - \lfloor \frac{n+r}{2} \rfloor-1} b(\lfloor \frac{n+r}{2} \rfloor
-1,k-1)=
C_{n - \lfloor \frac{n+r}{2} \rfloor-1} \sum_{i=1}^{\lfloor
\frac{n+r}{2} \rfloor-1} C_{i-1} b(\lfloor \frac{n+r}{2}
\rfloor-i-1,k-2),
$$
if $k\geq 2$, which holds by Lemma 4.3, thereby completing the
proof.
\Bx

\section{\large D-Wilf Classes and $\mathbf{h}$-Wilf Classes}

No discussion of restricted permutations is complete without
the discussion of {\it Wilf classes}, defined below, where
$S_n(T)$, $T \subseteq S_m$, is the set of $\pi \in S_n$
which avoid all patterns in $T$ and where $s_n(T)=|S_n(T)|$.

\begin{defn}
Let $S_1,S_2 \subseteq S_m$.  If $s_n(S_1)=s_n(S_2)$
for all $n \geq m$
then we say that $S_1$ and $S_2$ are in the same
Wilf class, or are Wilf equivalent.
\end{defn}

Since we have refined the investigation of restricted
permutations, we refine the notion of Wilf class with
the following definition, where
$S^k_n(T)$, $T \subseteq S_m$, is the set of $\pi \in S_n$
with exactly $k$ fixed points
which avoid all patterns in $T$ and where $s^k_n(T)=|S^k_n(T)|$.
We use $D_n(T)$ and $d_n(T)$ to represent 
$S_n^0(T)$ and $s_n^0(T)$, respectively.

\begin{defn}
Let $S_1,S_2 \subseteq S_m$.  If $d_n(S_1)=d_n(S_2)$ for all $n \geq m$
we say that $S_1$ and $S_2$ are in the same $D$-Wilf
class, or are $D$-Wilf equivalent.  If
$s_n^h(S_1) = s_n^h(S_2)$ for $h>0$,
for all $n \geq m,h$ we say
that $S_1$ and $S_2$ are in the same $h$-Wilf
class, or are $h$-Wilf equivalent.
\end{defn}

We will have need of the following lemma
in the proofs below.

\begin{lemma}  Let $\gamma \in S_n$ be given by $\gamma_i
=n+1-i$ for $1 \leq i \leq n$.  For $\pi \in S_n$, let $\pi^\star = \gamma
\pi
\gamma^{-1}$. Then, for all $\pi$, $\pi$ and $\pi^\star$ have the same
number of fixed points.  Furthermore, the number of occurrences
of the pattern $213$ (respectively $312$) in $\pi$ equals
the number of occurrences of the pattern $132$ (respectively $231$)
in $\pi^\star$.
\end{lemma}

\noindent
{\bf Proof.}  Since $\pi^\star$ is obtained from $\pi$ by
conjugation, $\pi$ and $\pi^\star$ have the same number of
fixed points.  If $i<j<k$ are such that $\pi_i \pi_j \pi_k$
is an occurrence of $213$ (resp. $312$), then $\gamma_k<\gamma_j<\gamma_i$
are such that $\pi^\star_{\gamma_k} \pi^\star_{\gamma_j}
\pi^\star_{\gamma_i}$ is an occurrence of $132$ (resp. $231$).
\Bx

We now state some results about refined Wilf classes.

\begin{thm}  There are exactly three $D$-Wilf classes
of patterns of length $3$.
\end{thm}

\noindent
{\bf Proof.}
Applying Theorem 5.1 with $r=0$ we see that,
in particular,
$d_n(321)=d_n(132)$.  Lemma 6.1 gives us
$d_n(132)=d_n(213)$ and $d_n(231)=d_n(312)$.

Lastly, we note that $d_4(123)=7$, $d_4(132)=6$,
and $d_4(231)=4$, thereby giving three $D$-Wilf
classes.
\Bx

\begin{thm}  For any $h>0$, there are exactly three $h$-Wilf classes
of patterns of length $3$.
\end{thm}

\noindent
{\bf Proof.}
We already have $s_n^h(132)=s_n^h(213)=s_n^h(321)$ 
and $s_n^h(231)=s_n^h(312)$ for all $1 \leq h \leq
n$.  We now show that for $h>0$, $\{s_n^h(132)\}_{n \geq 3}$,
$\{s_n^h(231)\}_{n \geq 3}$, and $\{s_n^h(123)\}_{n \geq 3}$ are different
sequences, thereby showing that there are three
$h$-Wilf classes for $h>0$.
 
First, since $s_5^1(123)=20$, $s_5^1(231)=16$,
$s_5^1(132)=13$, $s_5^2(123)=2$, $s_5^2(231)=8$, and
$s_5^2(132)=6$ we have exactly three $h$-Wilf classes
for $h=1,2$.   Next, 
note that for $h \geq 3$, $s_n^h(123)=0$
since with three fixed points we have a $123$ pattern.
We conclude by showing that
$s_{h+2}^h(132)=h+1$ and $s_{h+2}^h(312) \geq 2h+1$, thereby
giving exactly three $h$-Wilf classes for all
$h \geq 3$.

Let $\pi \in S_{h+2}^h(132)$.
To show that $s_{h+2}^h(132)=h+1$, note that we can only
have two entries of $\pi$ which are {\it not} fixed points.
Furthermore, $1$ cannot be a fixed point, for otherwise
all other entries must be fixed in order to avoid the
$132$ pattern.  Hence, we have the freedom to pick exactly
one of $2,3\dots,h+2$ to be a non-fixed point.  This gives
$h+1$ choices.

Next, we show that $s_{h+2}^h(312) \geq s_{h+1}^{h-1}(312)+2 $, which
gives $s_{h+2}^h(312) \geq 2h+1$ since
$s_3^1(312)=3$.  Consider the
following procedure.  Let $\pi \in S_{h+1}^{h-1}(312)$.
Let $\widehat{\pi}_i = \pi_i +1$ for $1 \leq i \leq h+1$.
Construct $X=\{ 1\widehat{\pi}_1 \widehat{\pi}_2 \cdots
 \widehat{\pi}_{h+1} \}
\cup
\{\pi_1\pi_2\cdots \pi_{h+1} (h+2): \pi_{1} \neq 1\}$.  
It is clear that these sets are disjoint and that 
$X \subseteq S_{h+2}^h(312)$.  All that remains to be shown
is that  $\mid \{\pi_1\pi_2\cdots \pi_{h+1} (h+2): \pi_{1} \neq 1\}
\mid = 2$.  To see this, note that $S_3^1(312) = \{132,321,213\}$
contains two elements with $\pi_{1} \neq 1$.  By construction
of the above procedure,
these two elements beget two elements in $S_4^2(312)$ such that $\pi_1
\neq 1$. This concludes the proof.
\Bx

We may further show, using Lemma 6.1, another result
concerning Wilf classes.  First, we remind
the reader of the following definition
from [Rob].

\begin{defn}
Let $S \subseteq S_m$ and let $T$
be a multiset of $S_m$.  Define $S_n(S;T)$ to be
the set of $\pi \in S_n$ which avoid all patterns
in $S$ and contain each element, including
multiplicities, in $T$ exactly once.
Let $s_n(S;T)=|S_n(S;T)|$.  
Let $S_1,S_2 \subset S_m$ and let $T_1,T_2$
be multisets of $S_m$.  If $s_n(S_1;T_1) = s_n(S_2;T_2)$
for all $n \geq m$ we say that
$(S_1;T_1)$ and $(S_2;T_2)$ are in the same almost-Wilf class,
or are almost-Wilf equivalent.  We drop the set
notation for a singleton set.
\end{defn} 

Next, we refine this definition.

\begin{defn}
Define $S^k_n(S;T)$ to be those permutations in $S_n(S;T)$
with exactly $k$ fixed points.  Let $s_n^k(S;T)
= |S^k_n(S;T)|$.  When $k=0$ we
write $D_n(S;T)$, and $d_n(S;T)$, respectively.
If $d_n(S_1;T_1) = d_n(S_2;T_2)$
for all $n \geq m$ we say that
$(S_1;T_1)$ and $(S_2;T_2)$ are in the same almost-$D$-Wilf class,
or are almost-$D$-Wilf equivalent. If $s^h_n(S_1;T_1) = s^h_n(S_2;T_2)$
for $h>0$ for all $n \geq m,h$ we
say that
$(S_1;T_1)$ and $(S_2;T_2)$ are in the same almost-$h$-Wilf class,
or are almost-$h$-Wilf equivalent.
\end{defn} 

\begin{thm}  Consider $D_n(\emptyset,\alpha)$, $\alpha \in S_3$.
For these permutations there are exactly four almost-$D$-Wilf classes.
\end{thm}

\noindent
{\bf Proof.}  We first prove that $d_n(\emptyset;321)=0$ for
all $n \geq 3$.  Let $cba$ be a $321$ pattern in $\pi \in S_n$.  
Write $\pi=WcXbYaZ$.  
In order to avoid another
$321$ pattern we see that for $w \in W$ and $x \in X$ we
must have $w,x<b$.  Furthurmore, for $y\in Y$ and $z\in Z$
we must have $y,z>b$.  Hence, $b$ is a fixed point.
Thus, the restriction of having exactly one $321$ pattern
implies a fixed point must be present.

Next, using Lemma 6.1,
we see that $d_n(\emptyset;132)=d_n(\emptyset;213)$
and $d_n(\emptyset;231)=d_n(\emptyset;312)$.
Lastly, since $d_5(\emptyset;123)=14, d_5(\emptyset;132)=8$,
and $d_5(\emptyset;231)=6$, we have exactly
four almost-$D$-Wilf classes.
\Bx

\section{\large Enumeration and Other Results}
\setcounter{equation}{0}

We start this section by tabulating $s_n^k(\alpha)$ for
$n \leq 8$ and $0 \leq k \leq 8$ for all $\alpha \in S_3$
using the fact that there are only
three $D$-Wilf classes. 

In the following tables let $n=0,1,\dots$ correspond to
row
$n$ and $k=0,1,2,\dots$ correspond to the $k^{\mathrm{th}}$ term in from
the left.  

\small
$$
\begin{array}{cccccccccccccccccc}
&&&&&&&&&1&&&&&\\
&&&&&&&&0&&1&&&&\\
&&&&&&&1&&0&&1&&&\\
&&&&&&2&&2&&0&&1&&\\
&&&&&6&&4&&3&&0&&1&\\
&&&&18&&13&&6&&4&&0&&1\\
&&&57&&40&&21&&8&&5&&0&&1\\
&&186&&130&&66&&30&&10&&6&&0&&1\\
&622&&432&&220&&96&&40&&12&&7&&0&&1\\
\end{array}
$$
\vskip 10pt
\centerline{\bf Values of 
$\mathbf{s_n^k(132)=s_n^k(321)=s_n^k(213)}$}

$$
\begin{array}{cccccccccccccccccc}
&&&&&&&&&1&&&&&\\
&&&&&&&&0&&1&&&&\\
&&&&&&&1&&0&&1&&&\\
&&&&&&1&&3&&0&&1&&\\
&&&&&4&&4&&5&&0&&1&\\
&&&&10&&16&&8&&7&&0&&1\\
&&&31&&44&&35&&12&&9&&0&&1\\
&&94&&146&&102&&59&&16&&11&&0&&1\\
&303&&464&&362&&180&&87&&20&&13&&0&&1\\
\end{array}
$$
\vskip 10pt
\centerline{\bf Values of 
$\mathbf{s_n^k(231)=s_n^k(312)}$}

$$
\begin{array}{cccccccccccccccccc}
&&&&&&&&&1&&&&&\\
&&&&&&&&0&&1&&&&\\
&&&&&&&1&&0&&1&&&\\
&&&&&&2&&3&&0&&0&&\\
&&&&&7&&4&&3&&0&&0&\\
&&&&20&&20&&2&&0&&0&&0\\
&&&66&&48&&18&&0&&0&&0&&0\\
&&218&&183&&28&&0&&0&&0&&0&&0\\
&725&&552&&153&&0&&0&&0&&0&&0&&0\\
\end{array}
$$
\vskip 10pt
\centerline{\bf Values of 
$\mathbf{s_n^k(123)}$}
\normalsize
From Corollary 3.3, Theorem 5.1,
and Lemma 6.1 we have $d_n(132)=d_n(321)=d_n(213)=F_n$
for $n \geq 1$, where $F_n$ is the $n^{\mathrm{th}}$ Fine number.

From the above triangles  it appears that
$d_n(231)<F_n$ and $d_n(123)>F_n$ for $n \geq 3$.  Unfortunately,
we were unable to prove the latter assertion.  However,
we can prove the former via a bijection similar to
one found in [Kra].

\begin{thm}
Let $F_n$ be the $n^{th}$ Fine number.
For all $n \geq 3$, $d_n(231) < F_n$. 
\end{thm}

\noindent
{\bf Proof.}  We first consider the
following bijection $\gamma:S_n(132)
\rightarrow SR_n$.  Let $\pi \in S_n(132)$,
$\pi=\pi_1\pi_2\cdots \pi_n$.  Let $s=
s_ns_{n-1}\cdots s_1 = \gamma(\pi)$
where $s_{i}$ is the number of entries
in $\pi_{i+1} \pi_{i+2} \cdots \pi_n$ which
are larger than $\pi_i$.  For example,
$\gamma(456312)=010012$.

The fact that $\pi$ is $132$-avoiding
guarantees that $s_{i-1} \leq s_{i}+1$.
Hence, $s \in SR_n$.  The inverse
bijection is obvious.

We now prove that if $s \not \in SR_n(0)$
then $\gamma^{-1}(s)$ is not a backward derangement.
By observing
that if $\pi \in S_n(132)$ is a backward derangement
then when $\pi$ is read from right to left it
becomes a member of $D_n(231)$, we can conclude
that $D_n(231) \leq |SR_n(0)| = F_n$.

We must first take care of the case where
$\gamma(\pi)$ produces $s_1=0$.
In this case, it is obvious that $\pi_1=n$
and hence $\pi$ is
not a backward derangement.  We now assume
that the first case of a double zero starts
with $s_{i+1}$, $i \neq 0$, so that
$s_{i+1}=s_i=0$.
The bijection implies that $\pi_i>\pi_{i+1}$
If there exists $\pi_k$, $k<i$, with $\pi_k<\pi_{i+1}$
then $\pi_k \pi_{i+1} \pi_i$ is
an occurrence of the pattern $132$, and hence is
not allowed.  Thus, all elements less than $\pi_{i+1}$ are
to the right of $\pi_{i+1}$.  This implies
that $\pi_{n+1-\pi_{i+1}}=\pi_{i+1}$.

We must now prove that the inequality is strict.
We do this by giving $\pi \in S_n(132)$ which
is not a backward derangement for which $\gamma(\pi)$
is nonsingular.
For $n$ odd
$\gamma(123\cdots n) = 012 \cdots n-1$
and for $n$ even
$\gamma(234\cdots (n-1)1n)=01123\cdots n-2$.
Noting that both of the similarity relations
produced by $\gamma$ are nonsingular and correspond
to permutations of $S_n(132)$ which are not
backward derangements completes the proof.
\Bx

We now turn our attention to the patterns $132, 321$, and $213$
(all of which are in the same $D$-Wilf class).
Using Theorem 5.1, we are able to derive some formulas
(recursive and ``closed" form) for $s_n^k(\alpha)$,
$\alpha \in \{132,321,213\}$.  First, we rederive
items 2 and 4 of Theorem 2.1.  

Let $r=n$ in Theorem 5.1.  Since $E_{n,n} = C_n$ we have,
using (4.1),
\begin{equation}
C_n = 2a_n(x) -xa_n(x) + (1-x)^2 a_{n-1}(x)
\end{equation}
Equating the constant terms gives us
$C_n=2d_n(\alpha)+d_{n-1}(\alpha)$ for $\alpha \in \{132,321,213\}$.
Hence, $C_n=2F_n + F_{n-1}$ is rederived.  From here the derivation
of item 4 in Theorem 2.1 uses either
a straightforward induction or telescoping sum.

We can, of course, use (7.1) to derive recurrences for
$s_n^k(\alpha)$, $\alpha \in \{132,321,213\}$, for $k \neq 0$. 
To this end, we have the following.

\begin{thm}  For $0 \leq k \leq n$, let $F^k_n = s_n^k(\alpha)$, 
$\alpha \in\{132,321,213\}$.  Define $F_n^{-1}=0$.  For $n \geq
2$, we have for
$0 \leq k \leq n$,
$$2F^k_n+F^k_{n-1} = F^{k-1}_n+2F^{k-1}_{n-1}-F^{k-2}_{n-1}.$$
\end{thm}

\noindent
{\bf Proof.}  Equate the coefficients
of $x^k$ in (7.1).
\Bx

We can now use a telescoping sum to show the following.

\begin{thm} Let $C_n$ and $F_n$ be the
$n^{th}$ Catalan and Fine number, respectively.
Let $\alpha \in \{132,321,213\}$.  For $n \geq 1$, 
$$s_n^1(\alpha) = \frac{1}{4} \sum_{i=0}^{n-1} \left( \frac{-1}{2}
\right)^{i}(C_{n-i}+3F_{n-i-1}).$$

\end{thm}

\noindent
{\bf Proof.}  Routine. \Bx

Other formulas for $s_n^k(\alpha)$, $\alpha \in \{132,321,213\}$,
$k \geq 2$, similar to the one in Theorem 7.3
(which are based on Catalan and Fine numbers) can be derived
using Theorem 7.2.  However, these formulas are rather
cumbersome.  Instead, we present some nicer formulas.  First, we
have the following lemma, which introduces a generalization of the Catalan
numbers.

\begin{lemma}  Let $n \geq 1$.
Write $a_n(x) = \sum_{k=0}^n b(n,k)(x-1)^k$.
Then
$$
b(n,k) = \frac{k+1}{n+1} {2n-k \choose n}.
$$
\end{lemma}

\noindent
{\bf Proof.}  Using Lemma 4.2 it is routine to
verify the stated formula.
\Bx

Since Lemma 7.4 gives us a natural generalization
of the Catalan numbers, we define
\begin{equation}
C_n^{(k)} = \frac{k+1}{n+1} {2n-k \choose n}.
\end{equation}

Using Lemma 7.4 it is easy to prove the following
formulas.

\begin{thm}
Let $\alpha \in \{132,321,213\}$.  For $n \geq 1$, $0 \leq k \leq n$, we
have 
\end{thm}
$$s^k_n(\alpha) = \sum_{j=0}^{n-k} (-1)^j {j+k \choose k}
C_n^{(k+j)},$$
i.e.
$$s^k_n(\alpha) =\sum_{j=0}^{n-k} (-1)^j \left(\frac{j+k+1}{n+1}\right)
  {2n-k-j \choose n} {j+k \choose k}.
$$

\noindent
{\bf Proof.}  We have $a_n(x) = \sum_{j=0}^n C_n^{(j)}(x-1)^j$.
Expanding $(x-1)^j$ we get
$$
a_n(x) = \sum_{j=0}^n \sum_{i=0}^j (-1)^{j-i} C_n^{(j)} {j \choose i} x^i.
$$
Equating the coefficients of $x^k$,
we get
$$
s^k_n(\alpha) = \sum_{j=k}^n (-1)^{j-k} C_n^{(j)} {j \choose k},
$$
which, after a change of variable, gives the desired result.
\Bx

We now harvest some other interesting results.

From Theorem 7.5 we have the following corollary which
relates $C_n^{(k)}$ and $C_n$.

\begin{cor} Let $C_n$ be the $n^{th}$ Catalan number
and let $C_n^{(k)}$ be as defined in (7.2).  Then, for $n \geq 0$,
$$
C_n=\sum_{k=0}^n \sum_{j=0}^{n-k} (-1)^j {j+k \choose k}
C_n^{(k+j)}
$$
\end{cor}

\noindent
{\bf Proof.}  Since $\sum_{k=0}^n s_n^k(\alpha)=C_n$
for any $\alpha \in S_3$, the result is immediate.
\Bx

We can also use Theorem 7.5 to rederive a formula
for the Fine numbers given in [Deu].

\begin{cor}  Let $F_n$ be the $n^{th}$ Fine number.
For $n \geq 1$,
$$
F_n = \sum_{j=0}^n (-1)^j C_n^{(j)}.
$$
\end{cor}

As a consequence of
$C_n=2F_n+F_{n-1}$ and the fact that
$F_n = d_n(\alpha)$ for $\alpha \in \{132,213,321\}$
we get the following theorem.

\begin{thm}  Let $\alpha \in \{132,213,321\}$ and
let $F(\pi)$ be the set of fixed points of
$\pi \in S_n(\alpha)$.  Define
$$
T_n(\alpha) = \{\pi \in S_n(\alpha) : F(\pi) \cap \{1,2,\dots, n-1\}
\neq \emptyset\}.
$$
Then $|T_n(\alpha)| = F_{n}$, where $F_{n}$ is the $n^{th}$
Fine number.
\end{thm}

\noindent
{\bf Proof.}  We have $C_n = s_n(\alpha)$ and
$F_n = d_n(\alpha)$.  If
$U_n(\alpha) = \{\pi_1 \cdots \pi_{n-1} n : \pi \in D_{n-1}(\alpha)\}$
then $S_n(\alpha)$ is the disjoint union
$S_n(\alpha) = D_n(\alpha) \cup T_n(\alpha) \cup U_n(\alpha)$.
Since $C_n=2F_n+F_{n-1}$, the result follows.
\Bx

\section*{\large References}
\footnotesize

[Deu] E. Deutsch, Dyck Path Enumeration,
{\it Discrete Math.} {\bf 204} (1999), 167-202.

[Fin] T. Fine, Extrapolation when Very Little is Known,
{\it Information and Control} {\bf 16} (1970),
331-359.

[Kra] C. Krattenthaler, Permutations with Restricted
Patterns and Dyck Paths, {\it Advances in Applied
Math.} {\bf 27} (2001), 510-530.

[Knu] D. Knuth, \underline{The Art of Computer Programming}, vol. 3,
Addison-Wesley, Reading, MA, 1973.

[Ric] D. Richards, Ballot Sequences and Restricted Permutations,
{\it Ars Combinatoria}, {\bf 25} (1988), 83-86.

[Rob] A. Robertson, Permutations Restricted by Two Distinct
Patterns of Length Three, {\it Advance in Applied Math.}
{\bf 27}, 548-561.

[Rog] D. G. Rogers, Similarity Relations on Finite Ordered Sets,
{\it Journal of Combinatorial Theory (A)} {\bf 23} (1977),
88-98.

[Sha] L. W. Shapiro, A Catalan Triangle, {\it Discrete Math.}
{\bf 14} (1976), 83-90.

[Slo] Sloane, N. J. A., The On-Line Encyclopedia of
Integer Sequences, {\tt http://www.research.att.com/
$\sim$njas/sequences}, A009766.

[SS] R. Simion and F. Schmidt, Restricted Permutations,
{\it European Journal of Combinatorics} {\bf 6} (1985), 383-406.

[SS2] F. Schmidt and R. Simion, Card
Shuffling and a Transformation in $S_n$, {\it Aeq. Math}
{\bf 44} (1992), 11-34.

[Str] V. Strehl, A Note on Similarity Relations,
{\it Discrete Math.} {\bf 19} (1977), 99-101.

[Wes] J. West, Permutations with Forbidden Subsequences and
Stack Sortable Permutations, Ph.D. Thesis, MIT, 1990.

\end{document}